\documentclass[a4paper,12pt]{article}
\usepackage{amsfonts}
\usepackage{amsmath}
\usepackage{amssymb}
\usepackage{graphicx}
\usepackage{epsf}
\usepackage{pgf,tikz}
\usepackage{epsfig}

\usepackage{color}
\usepackage{afterpage}
\usepackage{verbatim}
\newcommand{\norm}[1]{\left\Vert#1\right\Vert}

\newcommand{\D}{\mathrm{D}}

\usepackage{times}
\newtheorem{theo}{Theorem}

\begin{document}

\title{Computational design of acoustic materials using an adaptive
  optimization algorithm}

\author{L. Beilina  \thanks{
Department of Mathematical Sciences, Chalmers University of Technology and
Gothenburg University, SE-42196 Gothenburg, Sweden, e-mail: \texttt{\
larisa@chalmers.se}}
\and
E. Smolkin  \thanks{
%Department of Physics, Chalmers University of Technology and
%Gothenburg University, SE-42196 Gothenburg, Sweden, and
Department of Mathematics and Supercomputing, Penza State University, Penza, Russia, e-mail: \texttt{\  smolkin@chalmers.se}}
}
\date{}

\maketitle

\graphicspath{   % directories to search the graphics files
{FIGURES/}
{/nfs4/soleil.ce.chalmers.se/cab/ce/subsys/fs/backup/none/user/larisa/Acoustic/AdaptCloak2D/RESULTS/FIGURES/}
}

\begin{abstract}

We consider the problem of design of the acoustic structure of arbitrary geometry with prescribed desired properties.  We use
optimization approach for the solution of this problem and minimize
the Tikhonov functional on adaptively refined meshes. These meshes are
refined locally only in places where the acoustic structure should be
designed. Our special symmetric mesh refinement strategy together with
interpolation procedure allows the construction of the symmetric
acoustic material with prescribed properties.  Efficiency
of the
presented adaptive optimization algorithm is illustrated on the
construction of the symmetric acoustic material in two dimensions.

\end{abstract}

%\begin{keyword}

%nanophotonic  structure \sep  domain decomposition method \sep acoustic wave equation \sep
%    \sep adaptive finite element method \sep  finite difference method \sep  coefficient inverse problem \sep invisibility acoustic cloaking   \sep Tikhonov functional \sep Lagrangian approach
%\end{keyword}

\maketitle

\section{Introduction}

In this work we present a new adaptive optimization
algorithm which can construct  acoustic materials with
arbitrary geometry from desired scattering parameters.
We formulate our problem as a Coefficient Inverse Problem (CIP),
and our goal is
 to determine an unknown spatially distributed wave speed of
the acoustic wave equation from boundary measurements on the
adaptively refined meshes.
 To solve our CIP, we minimize the Tikhonov functional in order to
 find the wave speed distribution inside  designed domain which
 satisfies prescribed scattering properties.
  In the
case of numerical simulations of Section \ref{sec:numex} we formulate
these properties as obtaining as small as possible
reflections from the designed structure.  For minimization of the
Tikhonov functional we use Lagrangian approach and search for a
stationary point of it on the adaptively refined meshes. Compared with
other works on this subject \cite{B,BJ,btkm14} we need to refine
mesh locally only inside the known geometry.  For construction of a
new mesh we use symmetric mesh refinement strategy combined with the
interpolation procedure over the neighboring vertices for every
element in the mesh. This allows us finally to get acoustic material
of the symmetric structure.

To construct the desired   acoustic structure we formulate an adaptive
optimization algorithm which includes solution of the forward and
adjoint problems for the acoustic wave equation.  The domain
decomposition finite element/finite difference (FE/FD) method of
\cite{hybrid} is used for the computational solution of these
problems.  This method is implemented efficiently using
the software packages WavES \cite{waves} and PETSc
\cite{petsc}.  In the theoretical part of this work we present proof
of the energy estimate for a hyperbolic equation with one unknown
function - the wave speed- and different boundary conditions for the
case of our domain decomposition.  We illustrate efficiency of the
proposed method in numerical examples on the construction of new acoustic material in two
dimensions.  The goal of our numerical simulations is to reconstruct
the wave speed function of the hyperbolic equation from single
observations of the solution of this equation in space
and time which gives us as small  reflections as possible.
We note that the domain decomposition approach in this case is
particularly feasible for implementing of absorbing boundary
conditions \cite{EM}.

Developed in this work adaptive optimization method can be used in  construction
and design of new materials including nano-materials with so-called
cloaking properties, see \cite{cloak1,cloak2,cloak3}. To obtain
cloaking structures in all these works are used methods of
transformational optics which are based on the accordance between material
parameters and coordinate transformations. In the current work we
propose to use an adaptive optimization algorithm which is an
alternative approach for the construction of an approximate
cloaking. Depending on applications, this method can be used alone or
as a compliment to the method of transformational optics. Advantage of
a new technique compared to the transformational optics is fast
construction of any material of arbitrary geometry with desired
symmetric structure of any size.  This structure is not dependent on
the coordinate transformation and can be adapted to  desired
properties of the physical material.  The mesh size of the symmetric
structure can be defined as a parameter in the adaptive mesh
refinement procedure used in the optimization algorithm.  Thus, the
new algorithm allows efficiently compute a new material of any
symmetric structure with desired properties.  A first version of a
such algorithm was presented in \cite{AoA} for  design of a
nanophotonic structure.

The paper is organized as follows. In Section \ref{sec:modelhyb} we
present statements of the forward and inverse problems and in Section
\ref{sec:opt} we describe the Lagrangian approach for solution of our
CIP.  Stability estimates for the solution of forward and adjoint
problems are given in Section \ref{sec:energyerror1}.  In Section
\ref{sec:fem} we present the domain decomposition FEM/FDM to solve the
minimization problem of Section \ref{sec:opt}, and in Section
\ref{sec:ad_alg} we present an adaptive conjugate gradient algorithm
for the solution of our CIP.  Finally, in our concluding Section
\ref{sec:numex} we demonstrate  efficiency of the adaptive
optimization algorithm identifying the wave speed
function in two dimensions to construct material of  symmetric
structure which produce as small reflections as
possible.

\section{Statement of the forward and inverse problems}

\label{sec:modelhyb}

 Let $x = (x_1, x_2)$ denote a point in $\mathbb{R}^2$ in an
 unbounded domain $D$.
We model the wave propagation by the
following Cauchy problem for the scalar wave equation:
\begin{equation}\label{modelhyb}
%\begin{split}
\begin{cases}
\tilde{c}(x) \frac{\partial^2 u}{\partial t^2}  -  \triangle u   = 0 & ~ \mbox{in}~~\mathbb{R}^2 \times (0, \infty), \\
  u(x,0) = f_0(x), ~~~u_t(x,0) = 0 &~ \mbox{in}~~ D.
\end{cases}
%\end{split}
\end{equation}
Here, $u$ is the total wave pressure generated by the plane wave
$p(t)$ which is incident at $x_1 = x_0$ and propagates along $x_2$
axis,
$\tilde{c}(x)= \frac{1}{c(x)^2}$ is the isotropic function with
the spatially distributed wave speed $c(x)$.
%We note that in this work we use the single equation
%(\ref{modelhyb}) instead of the full Maxwell's equations, since in
%\cite{BMaxwell} was demonstrated numerically that in the similar
%numerical setting, as we will use in this note, other components of
%the electric field are negligible compared to the initialized one.  We
%also note that a scalar model of the wave equation was used
%successfully to validate reconstruction of the dielectric permittivity
%function with transmitted \cite{BK1, BK2} and backscattered
%experimental data \cite{ btkm14, btkm14b, KBKSNF, NBKF, NBKF2}.

\begin{figure}[tbp]
 \begin{center}
 \begin{tabular}{c}
{\includegraphics[scale=0.2, clip=]{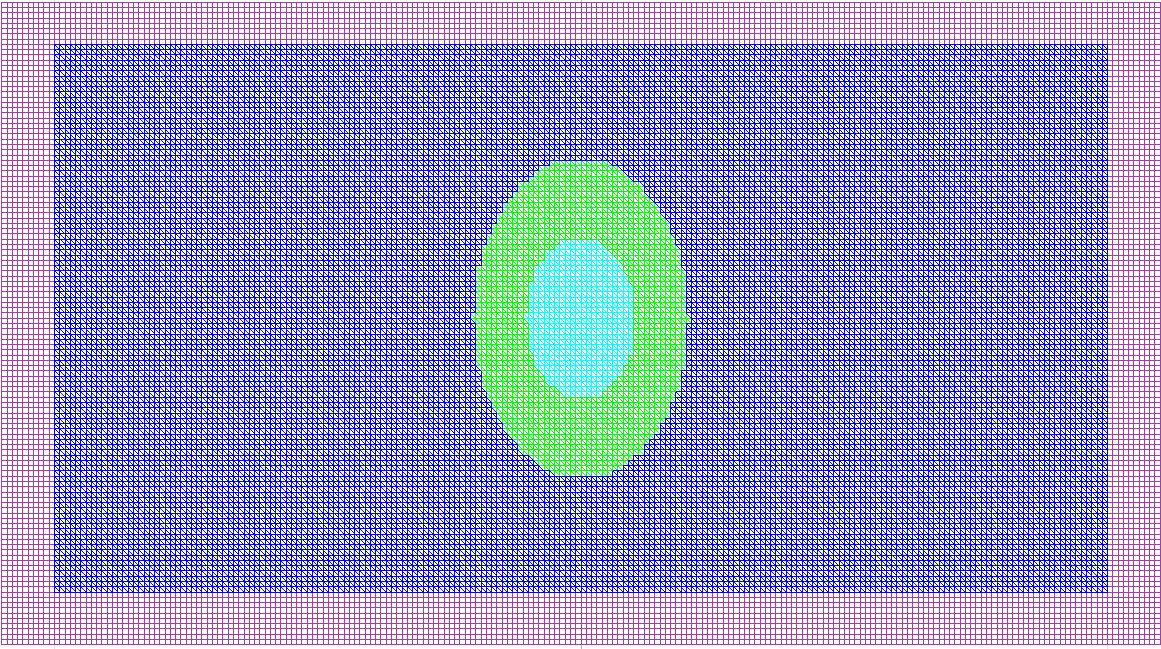}} \\
a) $D$ \\
%{\includegraphics[scale=0.19, clip=]{fdmmesh}} &
\begin{tikzpicture}[scale=0.7,x=0.74cm,y=0cm]
\node[anchor=north west,inner sep=0pt] at (0,0) {\includegraphics[scale=0.21, clip=]{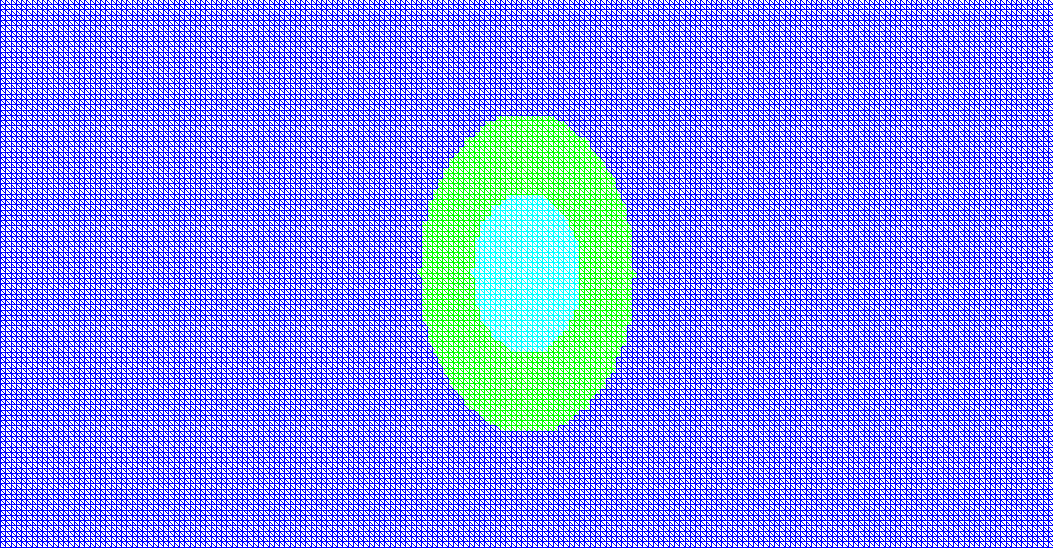}};
\draw (7.5,-15pt) node[anchor=north] {$G_2$};
\draw (7.5,-40pt) node[anchor= north] {$G_1$};
\draw (7.5,-65pt) node[anchor=north] {$G_0$};
\end{tikzpicture} \\
%{\includegraphics[scale=0.2, clip=]{FEMcoarsemesh}}\\
b)   $D_{FEM}$
 \end{tabular}
 \end{center}
 \caption{ \emph{a) Computational coarse FE/FD mesh  used in the domain
     decomposition of the domain $D = D_{FEM} \cup D_{FDM}$. b) The finite element
     mesh in $D_{FEM}$.}}
 \label{fig:0_1}
 \end{figure}

Let now $D \subset \mathbb{R}^{2}$ be a  bounded domain
with the boundary $\partial D$.  We use the notation
$D_T := D \times (0,T), \partial D_T := \partial D
\times (0,T), T > 0$ and assume that
\begin{equation}
 f_{0}\in H^{1}(D),  \tilde{c}(x) \in C^2(D).   \label{f1}
\end{equation}
For computational solution of (\ref{modelhyb}) we use the domain
decomposition finite element/finite difference (FE/FD) method of
\cite{hybrid} which was applied for the solution of different
coefficient inverse problems for the acoustic wave equation in
works \cite{B,hybrid,BJ,BCL}.  To apply method of \cite{hybrid} we
decompose $D$ into two regions $D_{FEM}$ and $D_{FDM}$ such that the
whole domain $D = D_{FEM} \cup D_{FDM}$,
% and $D_{FEM} \cap D_{FDM} =\emptyset$,
 see Figure \ref{fig:0_1}. In $D_{FEM}$ we use the finite
element method (FEM), and in $D_{FDM}$ we will use the Finite
Difference Method (FDM), see details in \cite{hybrid}.  Furthermore,
we decompose the domain $D_{FEM}$ into three regions $G_0, G_1, G_2$
such that $D_{FEM} = G_0 \cup G_1 \cup G_2$, where $G_0$ is the
innermost subdomain with the boundary $\partial G_0$, $G_1$ is the subdomain
where we want to design the acoustic material, and $G_2$ is the
outermost subdomain, see  Figure \ref{fig:0_1}-b).

 Let the boundary $\partial D$ be decomposed as $\partial D =\partial
 _{1} D \cup \partial _{2} D \cup \partial _{3} D$ where $\partial
 _{1} D$ and $\partial _{2} D$ are  top and bottom sides
 of the domain $D$, respectively, and $\partial _{3} D$ is the union of left and
 right sides of this domain.  At $S_T := (\partial_1 D \cup \partial_2
 D) \times (0,T)$ we have time-dependent observations.  We  define
 $S_{1} = \partial_1 D \times (0,T)$, $S_{1,1} = \partial_1 D \times
 (0,t_1]$, $S_{1,2} = \partial_1 D \times (t_1, T)$, $S_2 =
     \partial_2 D \times (0, T)$ and $S_3 = \partial_3 D \times (0,
     T)$, $S_4 = \partial G_0 \times (0, T)$.

We also introduce the following spaces of real valued  functions
\begin{equation}\label{spaces}
\begin{split}
H_u^1(D_T) &:= \{ w \in H^1(D_T):  w( \cdot , 0) = 0 \}, \\
H_{\lambda}^1(D_T) &:= \{ w \in  H^1(D_T):  w( \cdot , T) = 0\},\\
U^{1} &=H_{u}^{1}(D_T)\times H_{\lambda }^{1}(D_T)\times C\left( \overline{D}\right).
%\\
%U^{0} &=L_{2}\left(D_{T}\right) \times L_{2}\left(D_{T}\right) \times
%L_{2}\left( D \right).
\end{split}
\end{equation}

In our computations we have used the following model problem
\begin{equation}\label{model1}
%\begin{split}
\begin{cases}
\tilde{c} \frac{\partial^2 u}{\partial t^2}  -  \triangle  u   = 0&~ \mbox{in}~~ D_T, \\
u(x,0) = f_0(x), ~~~u_t(x,0) = 0 &~ \mbox{in}~~ D,     \\
\partial _{n}  u = p(x,t) & ~\mbox{on}~ S_{1,1}, \\
\partial _{n}  u =-\partial _{t} u & ~\mbox{on}~ S_{1,2}  \cup S_2, \\
%\partial _{n} u =-\partial _{t} u &~\mbox{on}~ S_2, \\
\partial _{n} u =0 &~\mbox{on}~ S_3 \cup S_4.
\end{cases}
%\end{split}
\end{equation}
In (\ref{model1}) we use the first order absorbing boundary conditions
\cite{EM}  and $p(x,t) \in L_2(S_{1,1})$.
We note that these conditions are exact in the case of
computations of Section \ref{sec:numex}, since in our computations we
initialize the plane wave orthogonal to the domain of propagation.

We choose the coefficient $\tilde{c}(x)$ in (\ref{model1})  such that
\begin{equation} \label{coefic}
\begin{cases}
\tilde{c} \left( x\right) \in \left [ 1, M \right], M=const. > 0,
& \text{ for }x\in  G_1, \\
 \tilde{c}(x) =1
& \text{ for }x\in  D_{FDM} \cup  G_2.
\end{cases}
\end{equation}

We consider the following inverse problem

\textbf{Inverse Problem (IP)}

\emph{Let the coefficient }$\tilde{c} \left( x\right)$\emph{\ in the
  problem (\ref{model1}) satisfies conditions (\ref{coefic}) and
  assume that }$ \tilde{c}\left( x\right) $\emph{\ is unknown in the
  domain }$G_1$\emph{. Determine the function }$ \tilde{c}\left(
x\right) $\emph{\ in (\ref{model1}) for }$x\in G_1$
\emph{\ assuming that the following function }$\widetilde u\left(
x,t\right) $\emph{\ is known}
\begin{equation}
  u\left( x,t\right) = \widetilde u \left( x,t\right), ~\forall \left( x,t\right) \in S_T.  \label{2.5}
\end{equation}

%The question of stability and uniqueness of \textbf{IP} is
%addressed in the recent work \cite{CristofolLiSoc}.

\section{Optimization method}

\label{sec:opt}

In this section we present the reconstruction method to solve inverse
problem \textbf{IP}. This method is based on the finding of the
stationary point of the following Tikhonov functional
\begin{equation}
F(u, \tilde{c}) = \frac{1}{2} \int_{S_T}(u - \widetilde{u})^2 z_{\delta }(t) dS dt +
\frac{1}{2} \gamma \int_{G_1}(\tilde{c}- \tilde{c}_0)^2~~ dx,
\label{functional}
\end{equation}
where $u$ satisfies the equations (\ref{model1}), $\tilde{c}_{0}$ is
the initial guess for $\tilde{c}$ (see details about choice of this
guess in  Section \ref{sec:numex} and \cite{btkm14,BOOK}),
$\widetilde{u}$ is the observed field at $S_T$, $\gamma > 0$ is the
regularization parameter and $z_{\delta }(t)$ is the compatibility
function in time and can be chosen as in \cite{btkm14}.

To find minimum of (\ref{functional}) we use the Lagrangian approach
\cite{B,BJ,btkm14} and define  the following Lagrangian in the week form
\begin{equation}\label{lagrangian1}
\begin{split}
L(v) &= F(u, \tilde{c})
-  \int_{D_T} \tilde{c} \frac{\partial
 \lambda }{\partial t} \frac{\partial u}{\partial t}  ~dxdt
+   \int_{D_T}( \nabla  u)( \nabla  \lambda)~dxdt  \\
& - \int_{S_{1,1}} \lambda p(x,t) ~dS dt    + \int_{S_{1,2} \cup S_2} \lambda \partial_t u ~dS dt,
%   + \int_{S_2} \lambda \partial_t u  ~dS dt, \\
\end{split}
\end{equation}
where $v=(u,\lambda, \tilde{c}) \in U^1$, and search for a stationary point
with respect to $v$ satisfying for all $\bar{v}= (\bar{u}, \bar{\lambda}, \bar{\tilde{c}}) \in U^1$
\begin{equation}
 L'(v; \bar{v}) = 0 ,  \label{scalar_lagr1}
\end{equation}
where $ L^\prime (v;\cdot )$ is the Jacobian of $L$ at $v$.

In order to find the Fr\'{e}chet derivative (\ref{scalar_lagr1}) of
the Lagrangian (\ref{lagrangian1}) we consider $L(v + \bar{v}) - L(v)~
\forall \bar{v} \in U^1$ and single out the linear part of the
obtained expression with respect to $\bar{v}$.  When we derive the
Fr\'{e}chet derivative we assume that in the Lagrangian
(\ref{lagrangian1}) functions in $v=(u,\lambda, \tilde{c}) \in U^1$ can be
varied independent on each others. We note that by doing so we get
the same Fr\'{e}chet derivative of the Lagrangian (\ref{lagrangian1})
as by assuming that functions $u$ and $\lambda$ are dependent on the
coefficient $\tilde{c}$, see details in Chapter 4 of \cite{BOOK}.
Similar to \cite{B,hybrid,BJ} we use conditions $\lambda \left(
x,T\right) =\partial _{t}\lambda \left( x,T\right) =0$ and imply such
conditions on the function $\lambda $ to deduce that $ L\left( u,\lambda,
\tilde{c} \right) :=L\left( v\right) =F\left( u, \tilde{c}\right).$ We
also use conditions (\ref{coefic}) on $\partial D$, together with
initial and boundary conditions of (\ref{model1}) to get that for all
$\bar{v} \in U^1$ we have
\begin{equation}
 L'(v; \bar{v}) = \frac{\partial L}{\partial \lambda}(v)(\bar{\lambda}) +  \frac{\partial L}{\partial u}(v)(\bar{u}) + \frac{\partial L}{\partial  \tilde{c}}(v)(\bar{\tilde{c}}) = 0,  \label{scalar_lagrang1}
\end{equation}
or
\begin{equation}\label{forward1}
\begin{split}
0 &= \frac{\partial L}{\partial \lambda}(v)(\bar{\lambda}) = \\
&- \int_{D_T} \tilde{c} \frac{\partial \bar{\lambda}}{\partial t} \frac{\partial u}{\partial t}~ dxdt
+  \int_{D_T}  (  \nabla u) (\nabla  \bar{\lambda}) ~ dxdt  \\
&- \int_{S_{1,1}} \bar{\lambda} p(x,t)~dS dt   \\
&+ \int_{S_{1,2} \cup S_2} \bar{\lambda} \partial_t u  ~dS dt,~~\forall \bar{\lambda} \in H_{\lambda}^1(D_T),
\end{split}
\end{equation}
\begin{equation} \label{control1}
\begin{split}
0 &= \frac{\partial L}{\partial u}(v)(\bar{u}) = \\
&\int_{S_T}(u - \widetilde{u})~ \bar{u}~ z_{\delta}~ d S dt- \int_{D}
\tilde{c} \frac{\partial{\lambda}}{\partial t}(x,0) \bar{u}(x,0) ~dx\\
&- \int_{S_{1,2} \cup S_2}
\frac{\partial{\lambda}}{\partial t} \bar{u} ~dS dt
 \\
&-  \int_{D_T} \tilde{c}  \frac{\partial \lambda}{\partial t} \frac{\partial \bar{u}}{\partial t}~ dxdt \\
&+ \int_{D_T} ( \nabla  \lambda) (\nabla  \bar{u})  ~ dxdt, ~\forall \bar{u} \in H_{u}^1(D_T),
\end{split}
\end{equation}
\begin{equation} \label{grad1new}
\begin{split}
0 &= \frac{\partial L}{\partial  \tilde{c}}(v)(\bar{\tilde{c}})
 =  -  \int_{D_T}  \frac{\partial \lambda}{\partial t} \frac{\partial u}{\partial t} \bar{\tilde{c}}~dxdt \\
&+\gamma \int_{G_1} (\tilde{c} - \tilde{c}_0) \bar{\tilde{c}}~dx,~ x \in D.
\end{split}
\end{equation}
We observe that (\ref{forward1}) is the weak formulation of the state equation
(\ref{model1}) and  (\ref{control1}) is the weak
formulation of the following adjoint problem
\begin{equation}\label{adjoint1}
\begin{cases}
\tilde{c} \frac{\partial^2 \lambda}{\partial t^2} -
  \triangle \lambda  = -  (u - \widetilde{u}) z_{\delta} &~  x \in S_T,   \\
\lambda(\cdot, T) =  \frac{\partial \lambda}{\partial t}(\cdot, T) = 0, \\
\partial _{n} \lambda = \partial _{t} \lambda & ~\mbox{on}~ S_{1,2} \cup S_2,
\\
%\partial _{n} \lambda =\partial _{t} \lambda & ~\mbox{on}~ S_2, \\
\partial _{n} \lambda =0 & ~\mbox{on}~ S_3 \cup S_4 \cup S_{1,1}.
\end{cases}
\end{equation}

 We define by   $u(\tilde{c}), \lambda(\tilde{c})$  exact solutions of the
forward and adjoint problems, respectively, for the known function $\tilde{c}$.
Then using the fact that
 exact solutions   $u(\tilde{c}), \lambda(\tilde{c})$  are sufficiently stable (see Chapter 5 of book \cite{lad} for details), we get from (\ref{lagrangian1})
\begin{equation}
F( u(\tilde{c}), \tilde{c}) = L(v(\tilde{c})),
\end{equation}
and the Fr\'{e}chet derivative of the Tikhonov functional can be written as
\begin{equation}\label{derfunc}
\begin{split}
F'(\tilde{c}) := &F'(u(\tilde{c}), \tilde{c})=  \frac{\partial F}{\partial \tilde{c}}(u(\tilde{c}), \tilde{c})
%(\bar{c})
  =  \frac{\partial L}{\partial \tilde{c}}(v(\tilde{c}))
%(\bar{c})
.
\end{split}
\end{equation}
Inserting (\ref{grad1new})  into  (\ref{derfunc}),
we get the following space-dependent function:
\begin{equation} \label{derfunc2}
\begin{split}
F'(\tilde{c})(x) &:= F'(u(\tilde{c}),\tilde{c})(x) =\\
&- \int_0^T \frac{\partial \lambda(\tilde{c})}{\partial t} \frac{\partial u(\tilde{c})}{\partial t} (x,t)~dt
+\gamma  (\tilde{c} - \tilde{c}_0)(x).
\end{split}
\end{equation}

\section{Stability estimates}

\label{sec:energyerror1}

The stability estimate for the forward problem (\ref{model1}) follows
from the stability estimate of \cite{hybrid} and can be derived using
the technique of \cite{lad}.  For analysis we first introduce the $L_2$ inner product and
the norm over $D_T$ and $D$, correspondingly, as
\begin{equation*}
\begin{split}
((a,b))_{D_T}    &= \int_{D} \int_0^T a b ~ dx dt,~ \|a \|_{L_2(D_T)}^2  = ((a,a))_{D_T}, \\
(a,b)_{D}    &= \int_{D} a b  ~dx,~ \| a \|_{L_2(D)}^2  = (a,a)_{D}.
\end{split}
\end{equation*}

 \textbf{Theorem}

 \emph{Assume that the condition (\ref{coefic}) for the function
   $\tilde{c}(x)$ holds. Let }$D \subset
 \mathbb{R}^{n}, n=2,3,$\emph{\ be a bounded domain with a piecewise smooth
   boundary }$\partial D$. \emph{\ For any }$ t\in \left(
 0,T\right) $\emph{\ we define }$D_{t}= \partial_1 D\times \left( 0,t_1\right)
 .$ \emph{\ Assume that there exists a solution }$u$\emph{\ of the problem (\ref{model1}).  Then $u \in
 H^{1}(D_{T})$
   is unique and there exists a positive constant}
 $A=A(\| \tilde{c} \|_{D}, t)$ \emph{\ such that the
   following energy estimate is true for all $t \in (0,T ]$}

\begin{equation}\label{estimate1}
\begin{split}
 \left \Vert \sqrt{\tilde{c}}~ \partial _{t} u(x,t)
 \right\Vert_{L_{2}\left( D \right) }^{2}
&+  \left \Vert \nabla u (x,t) \right\Vert _{L_{2}\left( D \right) }^{2} \\
 & \leq A \left[ \left\Vert p(x,t) \right\Vert _{L_{2}\left( D _{t}\right)}^{2}+
\left\Vert  \nabla f_{0}\right\Vert _{L_{2}\left( D\right) }^{2}\right] .
\end{split}
\end{equation}

\textbf{Proof.}

A proof of  this theorem follows from the stability estimate given in \cite{hybrid}.
$\square$

The stability result for the adjoint problem is obtained similarly as
for the forward problem, the only  difference is in the integration in
time $(t, T)$.
%We note also that we solve the adjoint problem backwards
%in time.

 \textbf{Theorem}

 \emph{Assume that the condition (\ref{coefic}) for the function
   $\tilde{c}(x)$ holds. Let }$D \subset
 \mathbb{R}^{n}, n=2,3$\emph{\ be a bounded domain with a piecewise smooth
   boundary }$\partial D$. \emph{\ For any }$ t\in \left(
 0, T\right) $\emph{\ we define by }$D_{t_a}=  (\partial_1 D \cup \partial_2 D) \times \left( t, T\right)
 .$ \emph{\ Assume that there exists a solution }$ \lambda $\emph{\ of the problem (\ref{adjoint1}) and a solution }$u $\emph{\ of the problem (\ref{model1}).  Then $\lambda \in
 H^{1}(D_{T})$
   is unique and there exists a positive constant}
 $B=B(\| \tilde{c} \|_{D}, t)$ \emph{\ such that the
   following energy estimate is true  for all $t \in (0,T ]$}

\begin{equation}\label{estimate2}
\begin{split}
 \left \Vert \sqrt{\tilde{c}}~ \partial _{t} \lambda(x,t)
 \right\Vert_{L_{2}\left( D \right) }^{2}
+  \left \Vert \nabla \lambda (x,t) \right\Vert _{L_{2}\left( D \right) }^{2}
 \leq B \left\Vert ( u - \tilde{u}) z_{\delta } \right\Vert _{L_{2}\left( D _{t_a}\right)}^{2}.
\end{split}
\end{equation}

\textbf{Proof.}

We multiply the equation in (\ref{adjoint1}) by $2 \partial_t \lambda$ and integrate
over $ D\times \left( t, T\right)$ to get
\begin{equation}\label{eq1_mod4adj}
\begin{split}
&\int \limits_{t}^{T}\int\limits_{D} 2~\tilde{c}~
\partial_{tt} \lambda ~\partial_t \lambda ~ dxd\tau
- \int \limits_{t}^{T}\int\limits_{D} 2
\nabla \cdot ( \nabla \lambda)~\partial_t \lambda~ dxd\tau \\
&=
- 2\int\limits_{t}^{T}\int\limits_{\partial_1 D \cup \partial_2 D}  (u - \tilde{u}) z_{\delta }  ~\partial_t \lambda ~dSd\tau.
\end{split}
\end{equation}

Next, we
integrate by parts in time the first term of (\ref{eq1_mod4adj}) and
noting zero initial condition in (\ref{adjoint1}),  we have
\begin{equation}\label{eq4_time1adj}
  \begin{split}
\int\limits_{t}^{T} \int\limits_{D} \partial _{t}
\left( \tilde{c} \partial_t \lambda^{2} \right)  dxd\tau
= -\int\limits_{D}\left( \tilde{c}\partial_t \lambda^{2} \right) \left( x,t\right) dx.
\end{split}
  \end{equation}

Next, we integrate by parts in space the second term of
(\ref{eq1_mod4adj}).  From (\ref{coefic}) it follows that $\tilde{c}=1$ on
$\partial D$. Thus, using (\ref{coefic}) and absorbing boundary
condition in (\ref{adjoint1}), we get
\begin{equation}\label{grad1adj}
\begin{split}
&2\int\limits_{t}^{T}\int\limits_{D} \nabla \cdot \left(\nabla  \lambda \right) \partial_t \lambda dx d\tau
=
  2\int\limits_{t}^{T}\int\limits_{\partial D}\left( \partial_t \lambda
  \right) \partial_n \lambda dS d\tau \\
  &-
  2  \int\limits_{t}^{T}\int\limits_{D}\left( \nabla \lambda \right) \left(
  \nabla \partial_t \lambda \right) dxd\tau \\
  &=
  2 \int\limits_{t}^{T}\int\limits_{\partial_1 D \cup \partial_2 D}\left( \partial_t \lambda
  \right)^2 ~dS d\tau - \int\limits_{t}^{T}\int\limits_{D} \partial _{t}|\nabla \lambda |^{2}dx d\tau.
\end{split}
\end{equation}

 Integrating last term of (\ref{grad1adj}) in time and using initial
 conditions of the equation (\ref{adjoint1}), we obtain
\begin{equation}\label{grad1_1adj}
\begin{split}
&\int\limits_{t}^{T}\int\limits_{D}  \partial _{t} |\nabla  \lambda| ^{2}dxd\tau
=\int\limits_{D} | \nabla  \lambda| ^{2}\left( x,T\right) dx
 - \int\limits_{D}  |\nabla \lambda |^{2}\left( x, t\right) dx \\
&= -\int\limits_{D}   | \nabla  \lambda|^{2}\left( x,t\right) dx.
\end{split}
\end{equation}

We
insert
(\ref{eq4_time1adj})-(\ref{grad1_1adj}) in (\ref{eq1_mod4adj}) to get
\begin{equation}\label{eq_mainadj1}
\begin{split}
&-\int\limits_{D}\left( \tilde{c} \partial_t \lambda^{2} \right) \left(
x,t\right) dx - \int\limits_{D} | \nabla \lambda |^{2}\left( x,t\right)
dx \\
&= 2 \left (\int\limits_{t}^{T}\int\limits_{\partial_1 D \cup \partial_2 D}\left(\partial_t
\lambda \right)^2  -
(u - \tilde{u} )z_{\delta } ~\partial_t \lambda \right)~~dS d\tau.
\end{split}
\end{equation}
The equation  above can be rewritten as
\begin{equation}\label{eq_mainadj2}
\begin{split}
&\int\limits_{D}\left( \tilde{c} \partial_t \lambda^{2} \right) \left(
x,t\right) dx + \int\limits_{D} | \nabla \lambda |^{2}\left( x,t\right)
dx \\
&= 2 \left (\int\limits_{t}^{T}\int\limits_{\partial_1 D \cup \partial_2 D}
( u -\tilde{u})z_{\delta } ~\partial_t \lambda
- \left(\partial_t
\lambda \right)^2\right )~~dS d\tau.
\end{split}
\end{equation}

%Using that $ 2\int\limits_{t}^{T}\int\limits_{\partial_1 D \cup \partial_2 D}\l%eft(
%\partial_t \lambda \right)^2 ~dS d\tau \geq 0$,  we
% get
%\begin{equation}\label{eq_mainadj}
%\begin{split}
%&\int\limits_{D}\left( \tilde{c} \partial_t \lambda^{2} \right) \left( x,t\right) dx
%+\int\limits_{D} | \nabla  \lambda |^{2}\left( x,t\right) dx \\
%&\leq 2\int\limits_{t}^{T}\int\limits_{\partial_1 D \cup \partial_2 D}
%|( u - \tilde{u} )z_{\delta }  |~|\partial_t \lambda|~dSd\tau.
%\end{split}
%\end{equation}

%Finally, to estimate the first term in the right hand side of
%(\ref{eq_main}) we use the arithmetic-geometric mean inequality $2ab
%\leq a^2 + b^2$ to obtain

Young's inequality applied to \eqref{eq_mainadj2}
directly leads to
\begin{equation}\label{mod4_5adj}
\begin{split}
&\int\limits_{D}
  \left(\tilde{c} \partial_t \lambda ^{2} + \left| \nabla \lambda \right|^{2} \right)(x,t)~ dx \\
  &\leq
B  \int\limits_{t}^{T}\int\limits_{\partial_1 D \cup \partial_2 D} |(\tilde{u} - u) z_{\delta }|^{2}(x, \tau)~dSd\tau,
\end{split}
\end{equation}
with a constant $B=0.5$
which is the desired result.

$\square$

\section{The finite element method in $D_{FEM}$}

\label{sec:fem}

As was mentioned above for the numerical solution of (\ref{model1}) we
use the domain decomposition FE/FD method of \cite{hybrid}.  Similarly
with this work, in our computations we decompose the finite
difference domain $D_{FDM}$ into squares, and the finite element domain
$D_{FEM}$ - into triangles. In $D_{FDM}$  we use the standard  finite
difference discretization of the equation (\ref{model1}) and obtain an
explicit scheme as in \cite{hybrid}.

For the finite element discretization of $D_{FEM}$ we define a
partition $K_{h}=\{K\}$ which consists of triangles. We define by $h$
the mesh function as $h|_{K}=h_{K}$, where $h_K$ is the local diameter of
the element $K$, and assume  the minimal
angle condition on the $K_{h}$ \cite{Brenner}. Let $J_{\tau }=\left\{ J\right\} $ be a partition of
the time interval $(0,\,T)$ into subintervals $J=(t_{k-1},\,t_{k}]$
of uniform length $\tau =t_{k}-t_{k-1}$.

To solve the state problem (\ref{model1}) and the
adjoint problem (\ref{adjoint1}) we define the finite element spaces,
$W_{h}^{u}\subset H_{u}^{1}\left( Q_{T}\right) $ and $W_{h}^{\lambda
}\subset H_{\lambda }^{1}\left( Q_{T}\right) $. First, we introduce the
finite element trial space $W_{h}^{u}$
\begin{equation*}
\begin{split}
W_{h}^{u} := &\{w\in H_{u}^{1}(Q_T):w|_{K\times J}\in P_{1}(K)\times
P_{1}(J), \\
&\forall K\in K_{h},~\forall J\in J_{\tau }\},
\end{split}
\end{equation*}
where $P_{1}(K)$ and $P_{1}(J)$ denote the set of linear functions on
$K$ and $J$, respectively. We also introduce the finite element test
space $W_{h}^{\lambda }$  as
\begin{equation*}
\begin{split}
W_{h}^{\lambda }:= &\{w\in H_{\lambda }^{1}(Q_T):w|_{K\times J}\in
P_{1}(K)\times P_{1}(J), \\
& \forall K\in K_{h},~\forall J\in J_{\tau }\}.
\end{split}
\end{equation*}
To approximate the function $\tilde{c}$, we use the space of
piecewise constant functions $C_{h}\subset L_{2}\left( D \right) $,
\begin{equation}
C_{h}:=\{u\in L_{2}(D ):u|_{K}\in P_{0}(K),~\forall K\in K_{h}\},
\notag
\end{equation}%
where $P_{0}(K)$ is the set of constant functions on $K$.

Setting $V_{h}=W_{h}^{u}\times W_{h}^{\lambda }\times C_{h}$, the
finite element method for  (\ref{scalar_lagr1}) now
reads: \emph{Find }$v_{h}\in V_{h}$\emph{, such that}
\begin{equation}
L^{\prime }(v_{h})(\bar{v})=0, ~\forall \bar{v}\in V_{h}.  \notag
\end{equation}

To find approximate solution $v_h \in V_{h}$ we need to solve the
forward problem \eqref{model1}, the adjoint problem \eqref{adjoint1}
and then find the discrete gradient $L^{\prime }_{\tilde{c}}(v_{h})$.  For
the fully discrete schemes of these equations we refer to
\cite{hybrid}.

\section{Adaptive conjugate gradient   algorithm}
\label{sec:ad_alg}

To compute minimum of the functional (\ref{functional}) we use the
adaptive conjugate gradient method (ACGM).  The regularization
parameter $\gamma$ in ACGM is computed iteratively via rules of
\cite{BKS}.  For the local mesh refinement we use a posteriori error
estimate of \cite {B,BJ} which means that the finite element mesh in
$D_{FEM}$ should be locally refined where the maximum norm of the
 Fr\'{e}chet derivative of the Lagrangian with respect to the
coefficient is large.  However, since our goal is to design material
inside the known domain $G_1$, we refine mesh only inside this domain.

Now we define
\begin{equation}\label{Bhm}
\begin{split}
  {g}^m(x) = - {\int_0}^T  \frac{\partial \lambda_h^m}{\partial t}
\frac{\partial E_h^m}{\partial t}~ dt   + \gamma^m (\tilde{c}_h^m - \tilde{c}_0),
\end{split}
\end{equation}
where $\tilde{c}_{h}^{m}$ is approximation of the function
$\tilde{c}_{h}$ on the iteration step $m$ in AGCM, $E_{h}\left(
x,t,\tilde{c}_{h}^{m}\right) ,\lambda _{h}\left( x,t,\tilde{c}_{h}^{m}
\right) $\ are computed by solving the state problem (\ref{model1}) and the
adjoint problem (\ref{adjoint1}), respectively, with
$\tilde{c}:=\tilde{c}_{h}^{m}$. In our computations of section \ref{sec:numex}
  we use the following algorithm.

\vspace{0.3cm}

\textbf{Algorithm  ~ (AGCM)}

\begin{itemize}
\item Step 0. Set number of mesh refinements $j:=0$. Choose initial
  mesh $K_{h}^j$ in $D_{FEM}$ and time partition $J_{\tau}^j$ of the time
  interval $\left( 0,T\right)$ as described in section
  \ref{sec:fem}. Start with the initial approximation
  $\tilde{c}_{h}^{0}= \tilde{c}_0$  at  $K_{h}^0$ and compute the sequences of
  $\tilde{c}_{h}^{m}$ via the following steps:

\item Step 1.  Compute solutions $E_{h}\left(
  x,t,\tilde{c}_{h}^{m}\right) $ and $\lambda _{h}\left(
  x,t,\tilde{c}_{h}^{m}\right) $ of state  (\ref{model1})
  and adjoint  (\ref{adjoint1})
problems, respectively, on $K_{h}^j$ and $J_{\tau}^j$.

\item Step 2.  Update the coefficient
  $\tilde{c}_h:=\tilde{c}_{h}^{m+1}$ on $K_{h}^j$ (only inside the discretized
  domain $G_1$) and $J_{\tau}^j$ using the conjugate gradient method
\begin{equation}\label{cgm}
\begin{split}
\tilde{c}_h^{m+1} &=  \tilde{c}_h^{m}  + \alpha^m d^m(x),
\end{split}
\end{equation}
where
\begin{equation*}
\begin{split}
 d^m(x)&=  -g^m(x)  + \beta^m  d^{m-1}(x),
\end{split}
\end{equation*}
with
\begin{equation*}
\begin{split}
 \beta^m &= \frac{\| g^m(x)\|^2}{\| g^{m-1}(x)\|^2},
\end{split}
\end{equation*}
where $d^0(x)= -g^0(x)$. In (\ref{cgm})
 the step size $\alpha$  in the gradient update is  computed as
\begin{equation}
\alpha^m = -\frac{((g^m, d^m)) }{\gamma^m {\norm{d^m}}^2},
\end{equation}
and the regularization parameter $\gamma^m$ at iteration $m$ is computed iteratively accordingly to \cite{BKS} as
\begin{equation}\label{iterreg}
\gamma^m = \frac{\gamma_0 }{ (m+1)^p},~~ p \in (0,1).
\end{equation}

\item Step 3. Stop computing $\tilde{c}_{h}^{m}$ and obtain the
  function $\tilde{c}_h$ at $M=m$ if either $\| g^{m}\|_{L_{2}( D_{FEM})}\leq
  \theta$ or norms $\|g^{m}\|_{L_{2}(D_{FEM})}$ are stabilized. Here
  $\theta$ is the tolerance in updates $m$ of gradient
  method. Otherwise set $m:=m+1$ and go to step 1.

\item Step 4.  Refine the mesh $K_h^j$ inside $G_1$ using symmetric mesh
  refinement procedure,  for example, as shown in Figure \ref{fig:7}.

\item Step 5. Set $j:=j+1$ and construct a new mesh $K_{h}^j$ in $D_{FEM}$ and a new partition $
J_{\tau}^j$ of the time interval $\left( 0,\,T\right)$ with the new time
step $\tau $ which should be chosen correspondingly to the CFL condition of \cite{CFL67}.

\item Step 6. Interpolate the approximation $\tilde{c}_h$ computed on
  the step 3, from every element $K^{j-1}$ on the previous space mesh
  $K_h^{j-1}$ to the new elements  $K^{j}$ in the mesh $K_h^{j}$, and obtain
  the initial guess $\tilde{c}_0$ on a new mesh.
  Set $m=1$ and return to step 1.

\item Step 7. Stop refinements of $K_h^j$ and $J_{\tau}^j$ if norms
  defined in step 3 either increase or stabilize, compared to the
  previous space mesh.

\end{itemize}

\section{Numerical Studies}
\label{sec:numex}

The goal of this section is to present possibility of the
computational design of an acoustic structure with the property to
generate as small reflections as possible.  This problem is equivalent
to \textbf{IP}. Thus, we will reconstruct a function $\tilde{c}(x)$
inside a domain $G_1$ using the ACGM algorithm of section
\ref{sec:ad_alg}.  We assume, that this function is known inside
$D_{FDM} \cup G_2$ and is set to be $\tilde{c}(x)=1$.

 Our computational geometry $D$
 is split into two geometries $D_{FEM}$ and $D_{FDM}$ as described in
 section \ref{sec:modelhyb}, see
 Figure \ref{fig:0_1}.
 We denote by $\partial D_{FEM}$  the outer boundary of  $D_{FEM}$
   and by $\partial D_{FDM}$ the inner boundary of $D_{FDM}$.
 We set the dimensionless computational domain
 $D$ as
 \begin{equation*}
 D = \left\{ x= (x_1,x_2) \in (-1.1, 1.1) \times (-0.62,0.62)\right\},
 \end{equation*}
and  the  domain $D_{FEM}$ as
 \begin{equation*}
 D_{FEM} = \left\{ x= (x_1,x_2) \in ((-1.0,1.0) \times (-0.52,0.52) \right\}.
 \end{equation*}
 The spatial mesh in $D_{FEM}$ and in $D_{FDM}$ consists of
 triangles and squares, respectively.
  We choose the initial mesh size $h=0.02$
 in $D = D_{FEM} \cup D_{FDM}$, as well as in the
 contiguous regions between FE/FD  domains.
We  also
decompose  the domain $\D_{FEM}$ into three different domains
$G_0, G_1, G_2$ such that $\D_{FEM} = G_0 \cup G_1 \cup G_2$ which are
intersecting only by their boundaries, see Figure \ref{fig:0_1}. The goal of our
numerical tests is to reconstruct the  function  $\tilde{c}$
of the  domain $G_{1}$ of Figure \ref{fig:0_1}
which produces as small reflections as possible.

  We initialize  a plane wave
  in
 $D$ in time $T=[0,2.0]$ such that
 \begin{equation}\label{f}
 \begin{split}
 p(t) =\left\{
 \begin{array}{ll}
 \sin \left( \omega t \right) ,\qquad &\text{ if }t\in \left( 0,\frac{2\pi }{\omega }
 \right) , \\
 0,&\text{ if } t>\frac{2\pi }{\omega }.
 \end{array}
 \right.
 \end{split}
 \end{equation}

 As for the forward problem   in $D_{FDM}$ we
 solve the problem (\ref{model1}) choosing $\tilde{c}=1$, and
in $D_{FEM}$ we  solve
 \begin{equation} \label{3D_1}
 \begin{split}
 \tilde{c} \frac{\partial^2 u}{\partial t^2} - \triangle u   &= 0~ \mbox{in}~~
 D_{{FEM} \times (0,T)},    \\
   u(x,0) = 0, ~~~u_t(x,0) &= 0~ \mbox{in}~~ D_{FEM},    \\
 u(x,t)|_{\partial D_{FEM}} &= u(x,t)|_{\partial D_{{FDM}_I}},\\
\partial_n u &= 0~ \mbox{on}~~ \partial G_0.
 \end{split}
 \end{equation}
Here, $\partial D_{{FDM}_I}$ denotes  internal
structured nodes of $D_{FDM}$ which have the same coordinates as  structured nodes
at the boundary $\partial D_{FEM}$, see details in \cite{hybrid}. We
note, that we use the boundary condition $\partial_n u = 0$ on
$\partial G_0$ which implies that waves are not penetrated into $G_0$.

We  also note that in $D_{FDM}$ the   adjoint problem  will be   the following wave equation  with
$\tilde{c}(x)=1$ for $ x  \in D_{FDM}$:
\begin{equation}\label{adjwaveeq}
\begin{split}
  \frac{\partial^2 \lambda}{\partial t^2} - \triangle \lambda   &=
  -  (u- \tilde{u}) z_{\delta}~~\mbox{in}~~ D_{FDM} \times (0,T),    \\
\lambda(x,T) = 0, ~~~\lambda_t(x,T) &=0~ \mbox{in}~~ D,     \\
\lambda(x,t)|_{\partial D_{FDM}} &= \lambda(x,t)|_{\partial D_{{FEM}_{I}}}, \\
\partial _{n} \lambda(x,t)& =0~ \mbox{on}~S_3 \cup S_{1,1}, \\
\partial _{n} \lambda(x,t)& = \partial_t \lambda~ \mbox{on}~S_{1,2} \cup S_{2},
%\nabla \cdot \lambda &= 0\text{ in }D ^{\prime }\subset D_{FDM}.
\end{split}
\end{equation}
which we solve using finite difference method.
 In $D_{FEM}$ we  solve the problem
\begin{equation} \label{adj3D_1}
\begin{split}
\tilde{c} \frac{\partial^2 \lambda}{\partial t^2} - \triangle \lambda &= 0~ \mbox{in}~~ D_{FEM} \times (0,T),    \\
  \lambda(x,T) = 0, ~~~\lambda_t(x,T) &= 0~ \mbox{in}~~ D_{FEM},    \\
\lambda(x,t)|_{\partial D_{FEM}} &= \lambda(x,t)|_{\partial D_{{FDM}_{I}}}, \\
\partial_n \lambda &= 0 ~ \mbox{on}~~  S_4, \\
\end{split}
\end{equation}
using finite element method.  Here, $\partial
D_{{FEM}_{I}}$ denotes internal structured nodes of
$D_{FEM}$ lying on
the inner boundary $\partial D_{FDM}$ of $D_{FDM}$, see details in \cite{hybrid} for the
exchange procedure between FE/FD solutions.

As initial guess $\tilde{c}_0(x)$ we take different constant values
of the function $\tilde{c}(x)$ inside domain of $G_1$ of Figure
\ref{fig:0_1} on the coarse non-refined mesh, and we take
$\tilde{c}(x)=1.0$ everywhere else in $D$. We choose three different
constant values of $\tilde{c}_0(x)= \{1.5,2.0,2.5\}$ inside $G_1$.
  We define that the minimal and maximal values of the
   function $\tilde{c}(x)$  belong
   to the following set $M_{\tilde{c}} $ of admissible parameters
 \begin{equation}\label{admpar}
 \begin{split}
  M_{\tilde{c}} := \left \{\tilde{c}\in C(\overline{D })|  1 \leq \tilde{c}(x)\leq  \max_{G_1} \tilde{c}_0(x) \right \}.
 \end{split}
 \end{equation}
 The time step is chosen to be $\tau=0.002$ which satisfies
the CFL condition \cite{CFL67}.

 \subsection{Reconstructions}

 We generate data at the observation points at $S_{T}$ by solving the
 forward problem (\ref{model1}) in the time interval $t=[0,2.0]$, with
 function $p(t)$ given by (\ref{f}) and for different values of
 $\omega= \{40,60,80,100\}$.  To generate non-reflected data
 $\tilde{u}$ at $S_{T}$ we take the function $\tilde{c}(x)=1$ for all
 $x$ in $D$ and solve the problem (\ref{model1}) with a plane wave
 (\ref{f}) and $\omega= \{40,60,80,100\}$.
We regularize the solution of the inverse problem by starting
computations with regularization parameter $\gamma=0.01$ in
(\ref{functional}) and then updating this parameter iteratively in
ACGM by formula (\ref{iterreg}).  Computing  the regularization
parameter in this way is optimal  for our problem.  We refer to
\cite{Engl} for different techniques for choice of a regularization
parameter.

\begin{figure}[tbp]
 %\begin{floatingfigure}[h]{3.5cm}
\begin{center}
\begin{tabular}{cc}
{\includegraphics[scale=0.375, clip=]{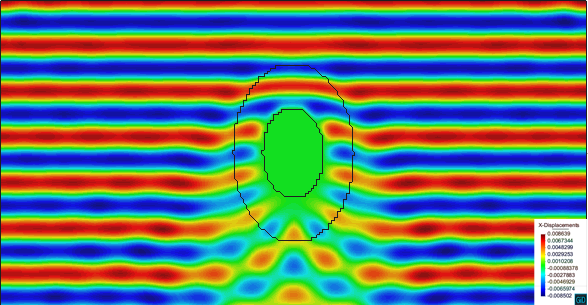}} &
{\includegraphics[scale=0.375, clip=]{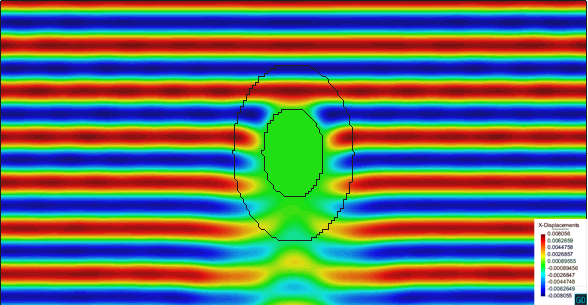}} \\
a)  $\omega= 40$ & b) $\omega= 40$\\
{\includegraphics[scale=0.375, clip=]{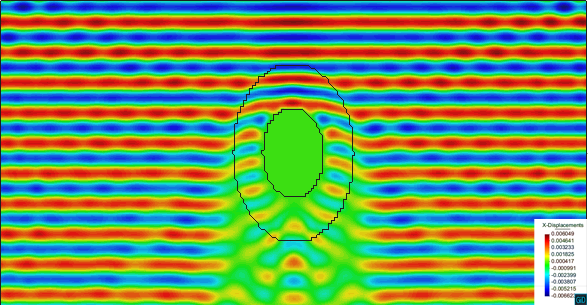}} &
{\includegraphics[scale=0.375, clip=]{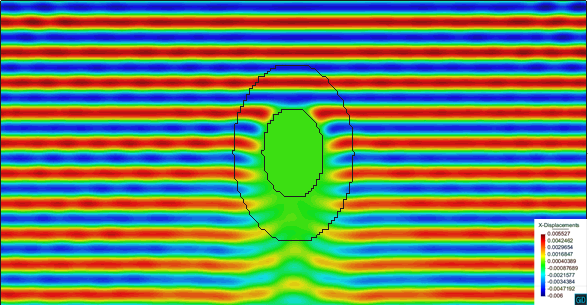}} \\
c) $\omega= 60$ & d)  $\omega= 60$ \\
{\includegraphics[scale=0.375, clip=]{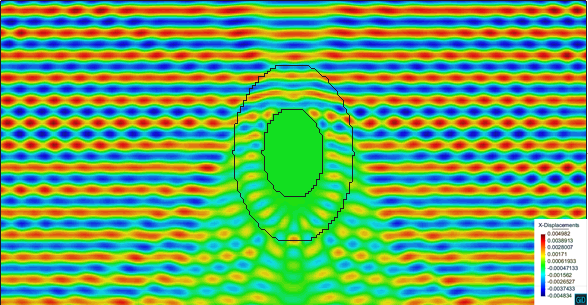}} &
{\includegraphics[scale=0.375, clip=]{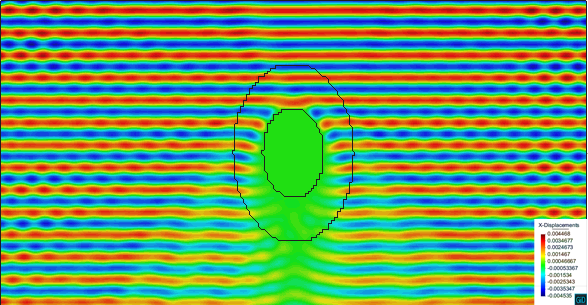}} \\
e) $\omega= 80$ & f) $\omega= 80$\\
{\includegraphics[scale=0.375, clip=]{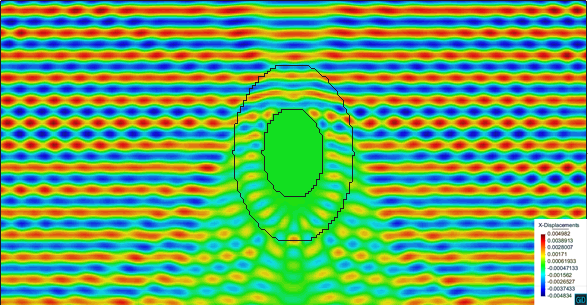}} &
{\includegraphics[scale=0.375, clip=]{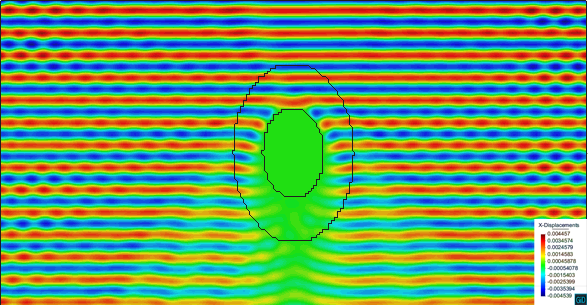}} \\
g) $\omega= 100$ & h) $\omega= 100$
\end{tabular}
\end{center}
\caption{{\protect\small \emph{Real part of the computational solution
      of (\ref{model1}) after the Fourier transform in time at different
      frequencies $\omega$: a),c),e),g) on the coarse mesh with
      $\tilde{c}_0=1.5$ in $G_1$; b),d),f),h) on the three times
      refined mesh with optimized $\tilde{c}$ in $G_1$. Optimized $\tilde{c}$
      for different frequencies $\omega$ is presented in Figures \ref{fig:4},
      \ref{fig:5}.}}}
\label{fig:2}
\end{figure}

Figure \ref{fig:2} shows real part of the Fourier transform of the
time-dependent solution $u(x,t)$ of \eqref{model1} when the initial
guess for $\tilde{c}$ was $\tilde{c}_0=1.5$ in all points of $G_1$
(left figures), and after application of the adaptive optimization
algorithm on three times refined mesh in $G_1$ (right  figures) for
different values of $\omega$ in (\ref{f}).  All right figures in
Figure \ref{fig:2} show significant reduction of backscattered
reflections for all tested frequencies compared with left figures.

\begin{figure*}[tbp]
\begin{center}
\begin{tabular}{cc}
%  {\includegraphics[scale=0.35, clip=]{Pictures/rec_1-1.png}}
{\includegraphics[scale=0.25, clip=]{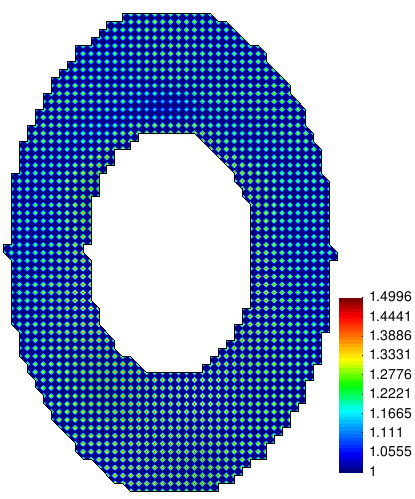}}
  &
%{\includegraphics[scale=0.25, clip=]{Pictures/rec_1-2.png}}
{\includegraphics[scale=0.125, clip=]{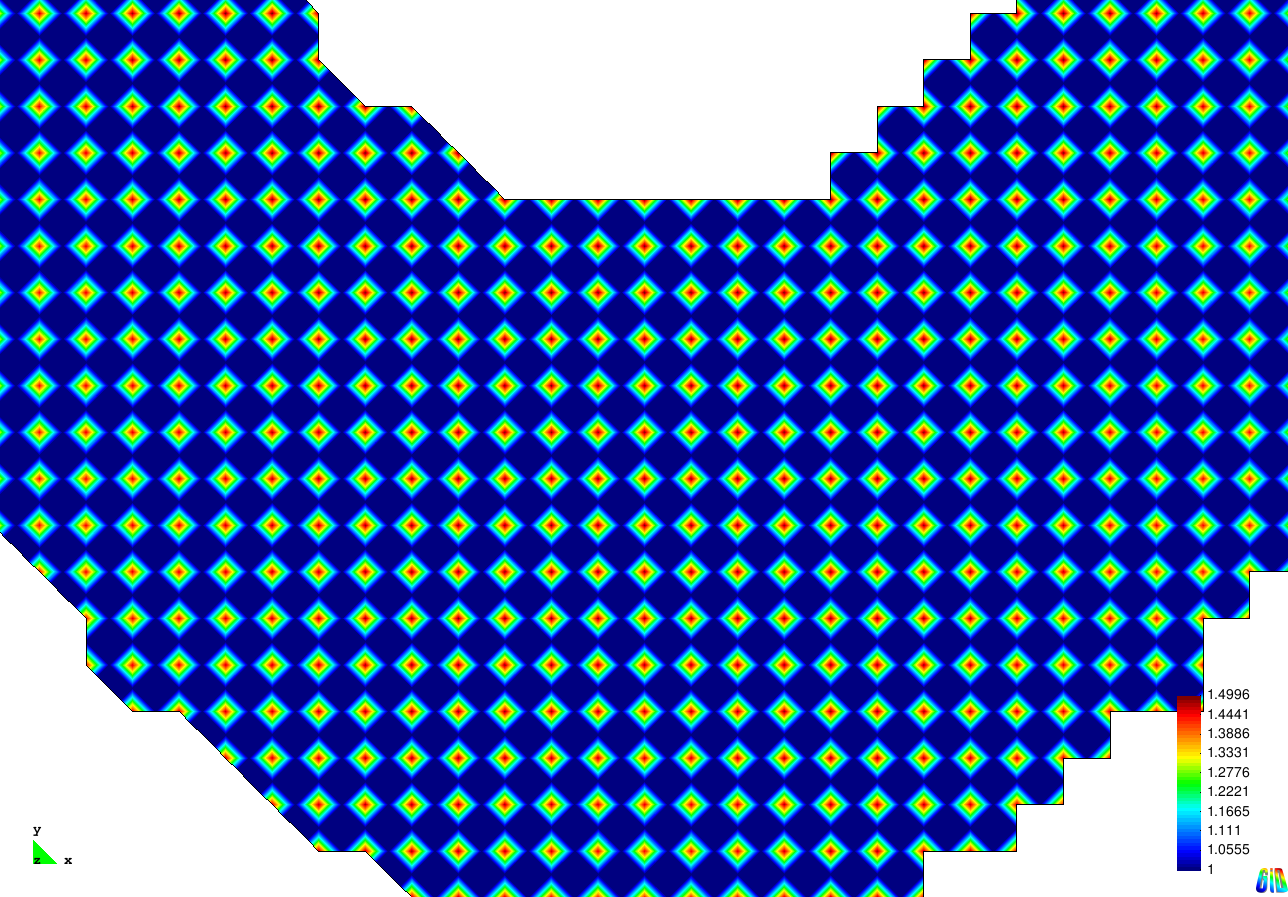}}
\\
{ a)  $\omega = 40, j=3$} & {zoomed}  \\
%{\includegraphics[scale=0.35, clip=]{Pictures/rec_2-1.png}}
{\includegraphics[scale=0.25, clip=]{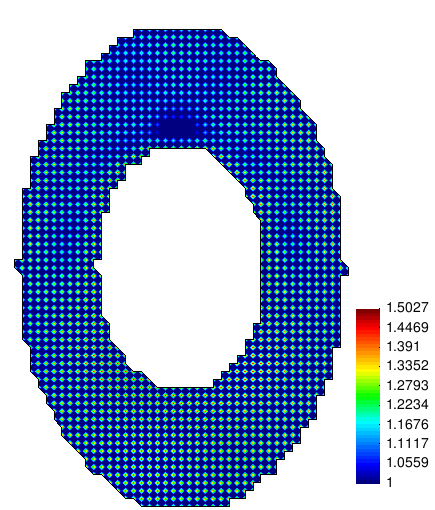}}
&
%{\includegraphics[scale=0.25, clip=]{Pictures/rec_2-2.png}}
{\includegraphics[scale=0.125, clip=]{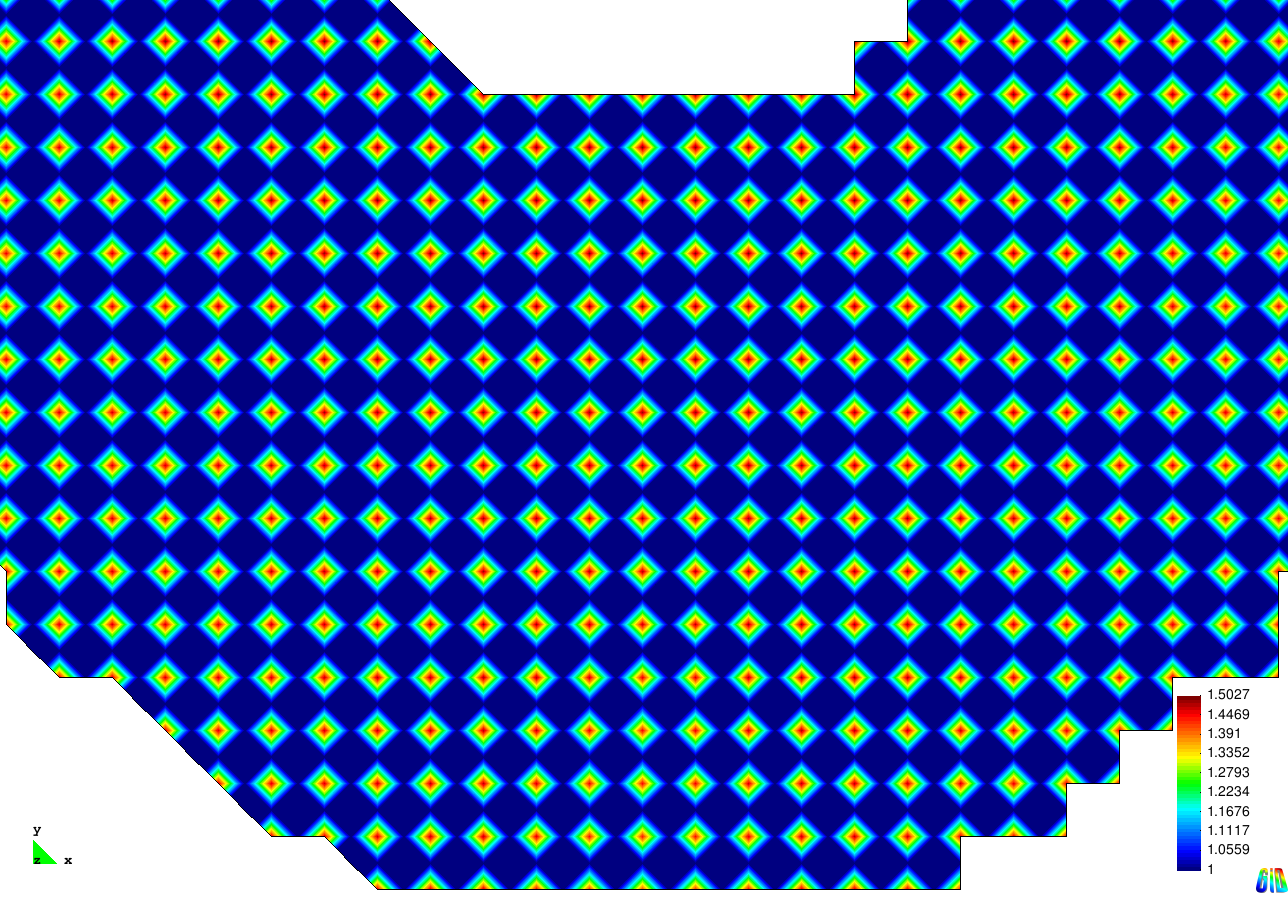}}
\\
{ b)    $\omega = 60, j=3$  } & {zoomed}
\end{tabular}
\end{center}
\caption{{\protect\small \emph{Reconstructed functions $\tilde{c}$ in
      $G_1$ on three times adaptively refined mesh $K_h^j,j=3,$ for
      different frequencies. Right figures present zoomed central bottom
      part of the domain $G_1$.}}}
\label{fig:4}
\end{figure*}

\begin{figure*}[tbp]
\begin{center}
\begin{tabular}{cc}
%  {\includegraphics[scale=0.35, clip=]{Pictures/rec_3-1.png}}
{\includegraphics[scale=0.25, clip=]{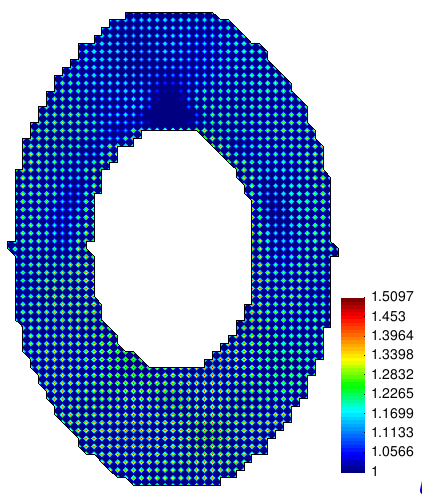}}
  &
%{\includegraphics[scale=0.25, clip=]{Pictures/rec_3-2.png}}
{\includegraphics[scale=0.125, clip=]{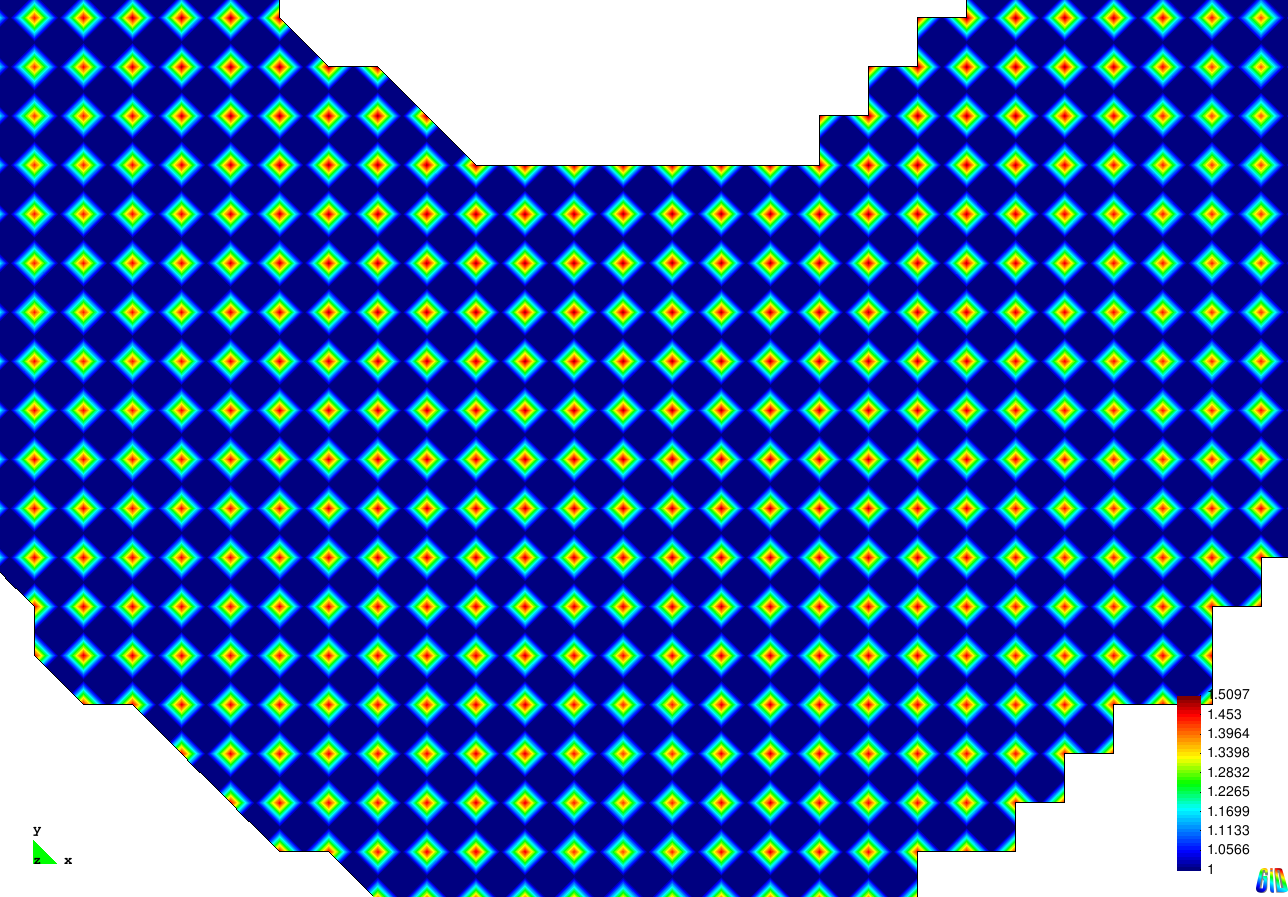}}
\\
{ a)    $\omega = 80, j=3$  } & {zoomed} \\
%{\includegraphics[scale=0.35, clip=]{Pictures/rec_4-1.png}}
{\includegraphics[scale=0.25, clip=]{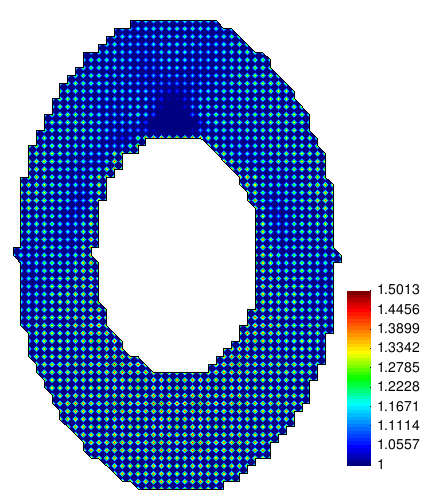}}
&
%{\includegraphics[scale=0.25, clip=]{Pictures/rec_4-2.png}}
{\includegraphics[scale=0.125, clip=]{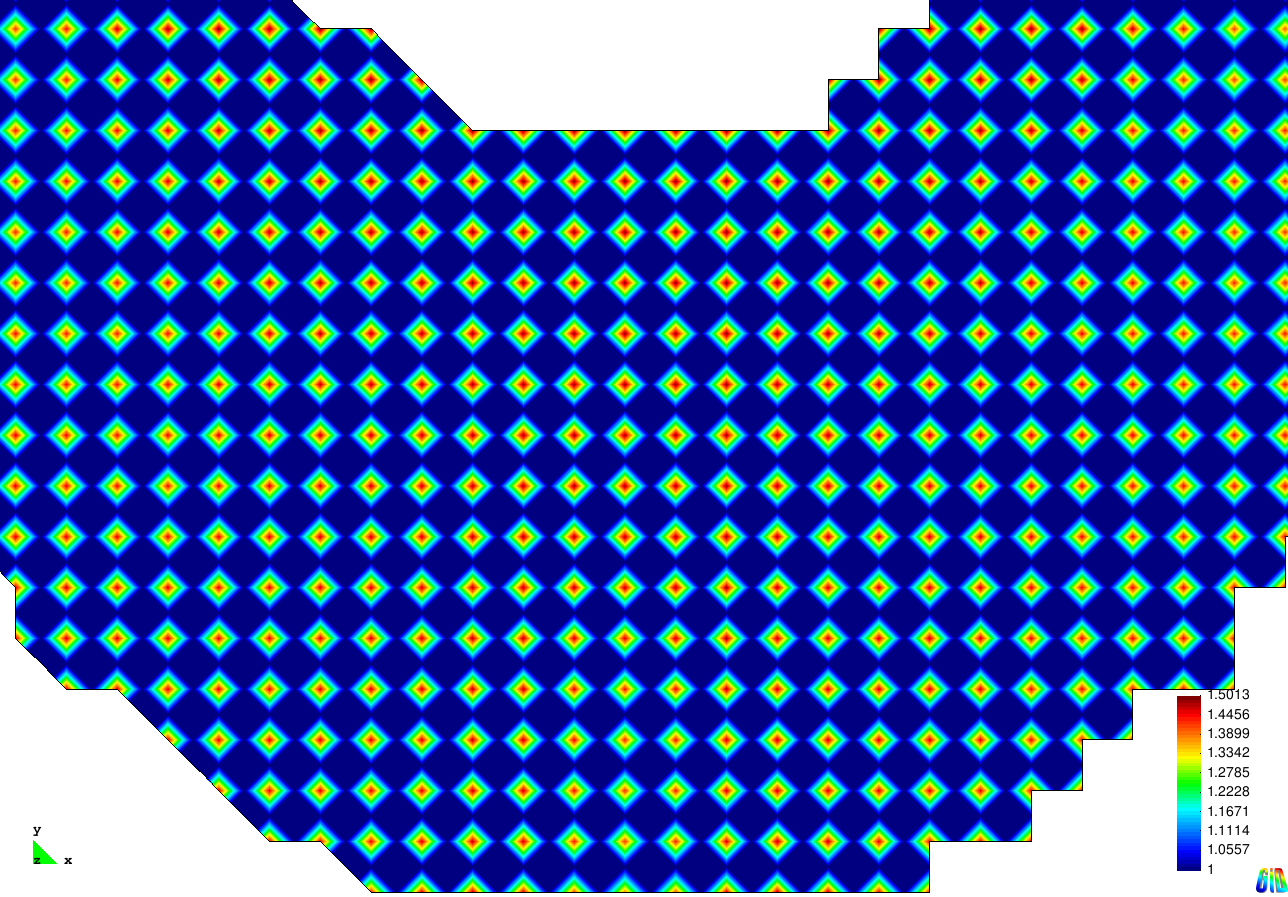}}
\\
{ b)    $\omega = 100, j=3$  } & {zoomed}
\end{tabular}
\end{center}
\caption{{\protect\small \emph{Reconstructed functions $\tilde{c}$ in
      $G_1$ on three times adaptively refined mesh $K_h^j, j=3,$ for
      different frequencies. Right figures present zoomed central bottom
      part of the domain $G_1$.}}}
\label{fig:5}
\end{figure*}

\begin{figure*}[tbp]
\begin{center}
\begin{tabular}{|c|c|c|c|}
\hline
{ Mesh} &{ Zoomed mesh} & {Reconstruction} & {Zoomed reconstruction}  \\
\hline
&&&\\
{\includegraphics[scale=0.2, clip=]{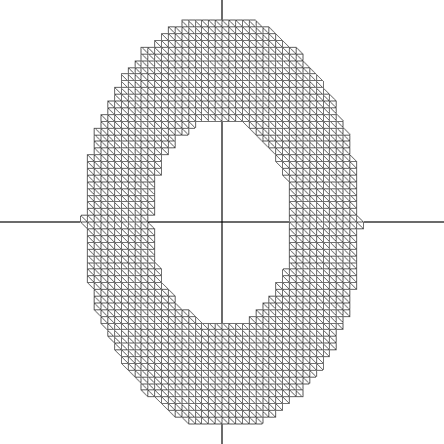}} &
{\includegraphics[scale=0.2, clip=]{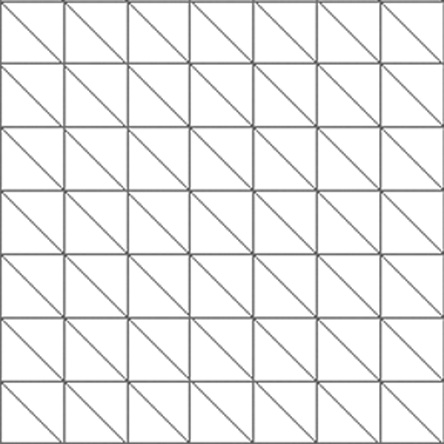}} &
{\includegraphics[scale=0.15, clip=]{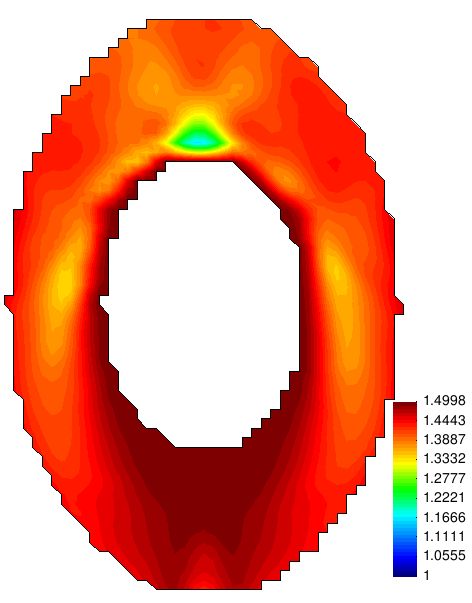}} &
{\includegraphics[scale=0.15, clip=]{coarseomega60.png}} \\
&&&\\
\hline
&&&\\
{\includegraphics[scale=0.2, clip=]{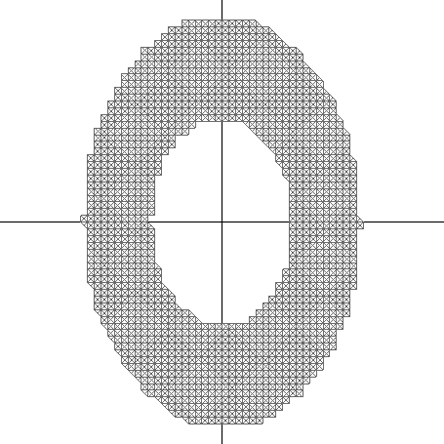}} &
{\includegraphics[scale=0.2, clip=]{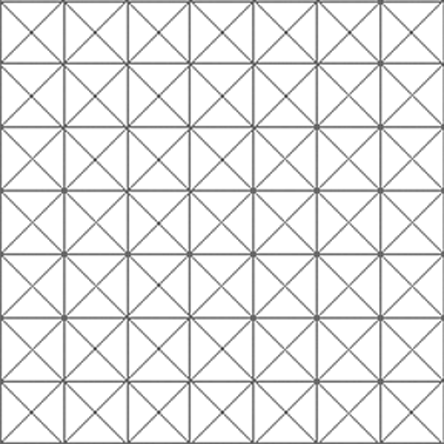}} &
{\includegraphics[scale=0.15, clip=]{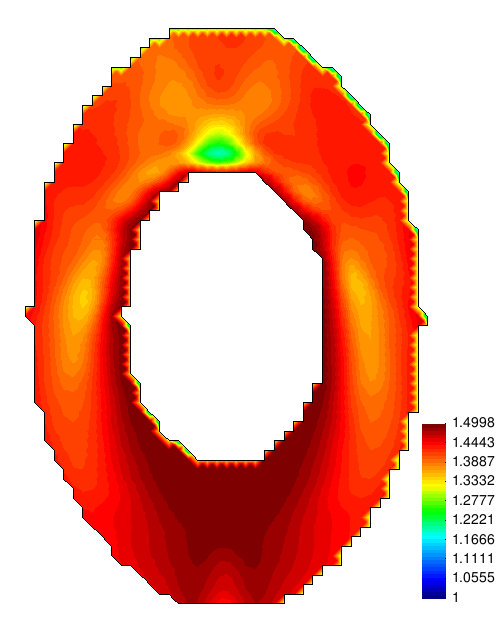}} &
{\includegraphics[scale=0.15, clip=]{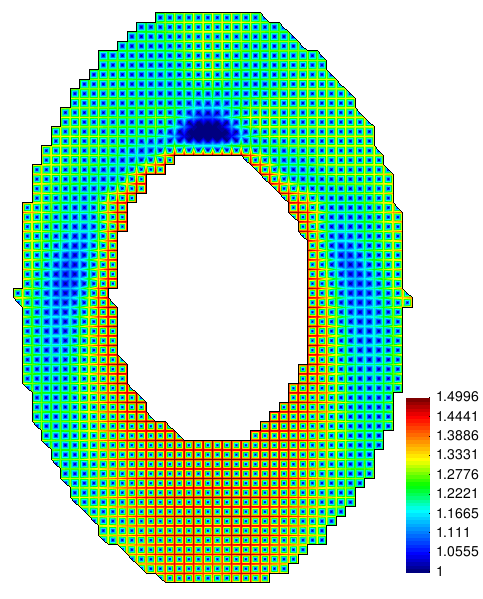}} \\
&&&\\
\hline
&&&\\
{\includegraphics[scale=0.2, clip=]{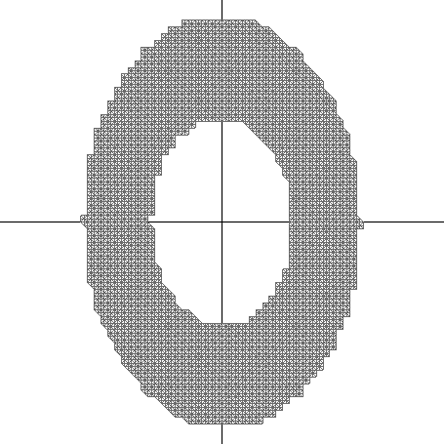}} &
{\includegraphics[scale=0.2, clip=]{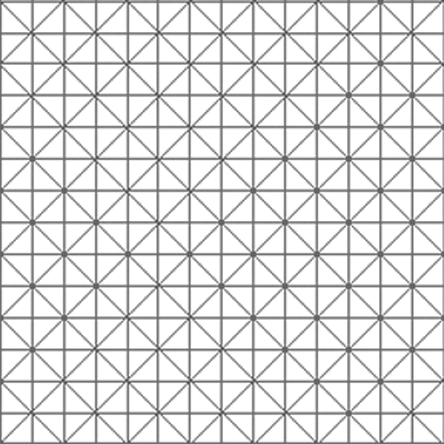}} &
{\includegraphics[scale=0.15, clip=]{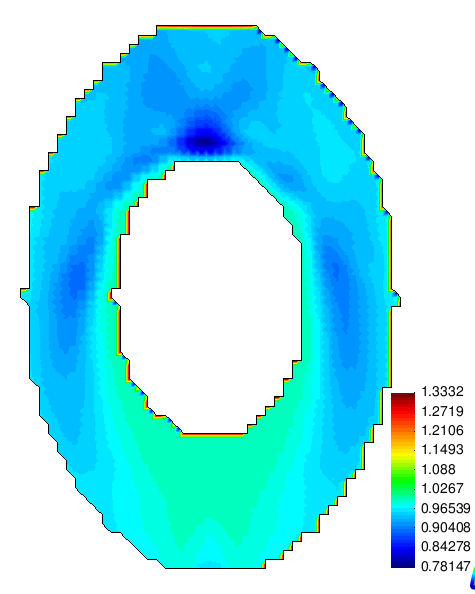}} &
{\includegraphics[scale=0.15, clip=]{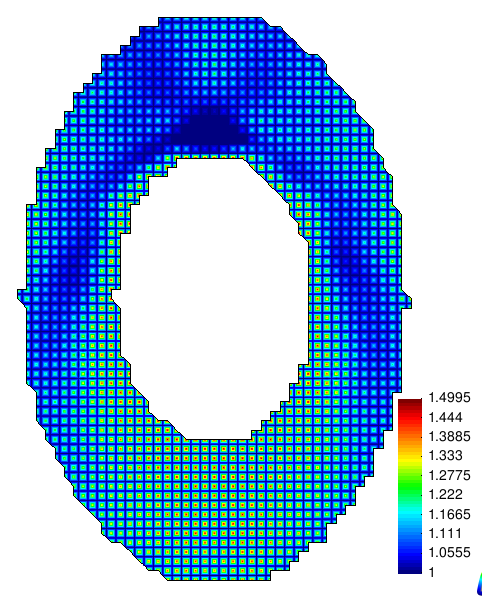}} \\
&&&\\
\hline
&&&\\
{\includegraphics[scale=0.2, clip=]{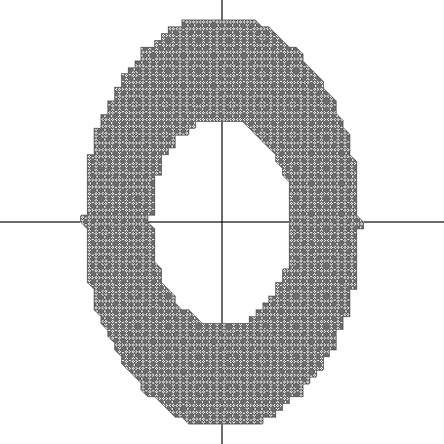}} &
{\includegraphics[scale=0.2, clip=]{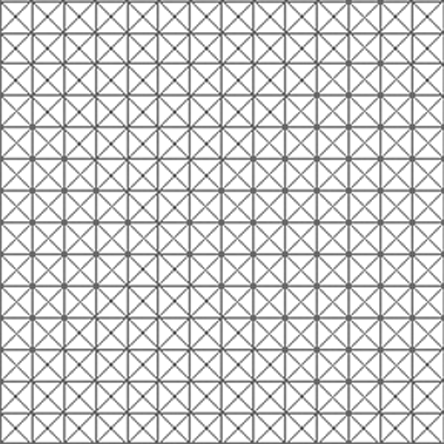}} &
{\includegraphics[scale=0.15, clip=]{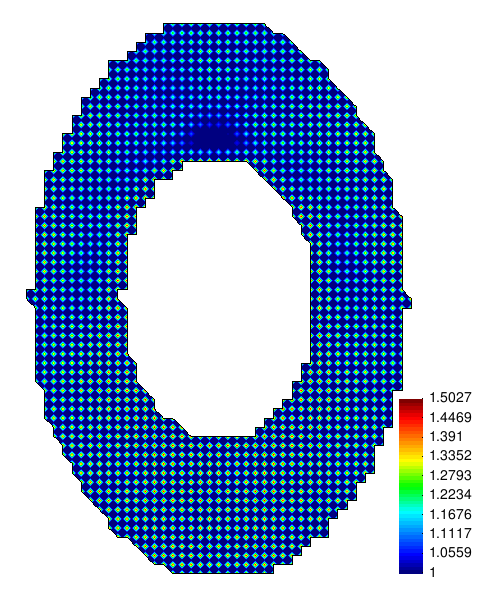}} &
{\includegraphics[scale=0.15, clip=]{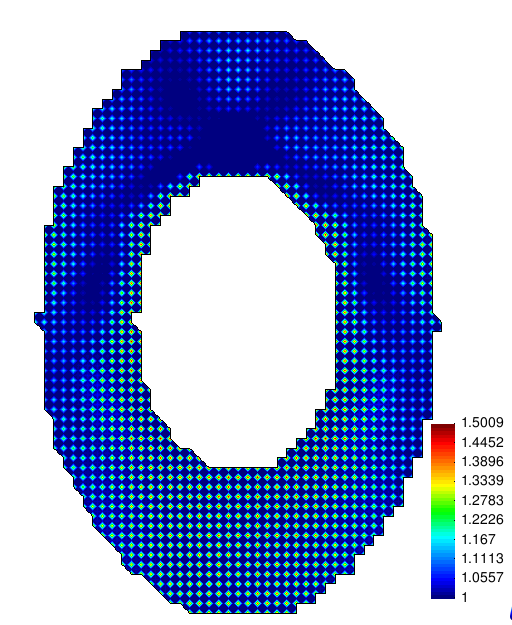}} \\
&&&\\
\hline
%&&&\\
%{\includegraphics[scale=0.2, clip=]{Pictures/res_4-1}} &
%{\includegraphics[scale=0.2, clip=]{Pictures/mesh_4}} &
%{\includegraphics[scale=0.2, clip=]{Pictures/res_4-2}} &
%{\includegraphics[scale=0.2, clip=]{Pictures/res_4-3}} \\
%&&&\\
%\hline
\end{tabular}
\end{center}
\caption{{\protect\small \emph{Reconstructed functions $\tilde{c}$ in
      $G_1$ for $ \omega=60$ in \eqref{f} on refined meshes $K_h^j,
      j=0,1,2,3$.  Left reconstructions: the optimized solution
      obtained on $K_h^0$ is interpolated on $K_h^j, j=1,2,3$.  Then
      the interpolated  $\tilde{c}$ on $K_h^3$ is taken as an initial guess and
      optimized further to get final reconstructed material shown on $K_h^3$.  Right
      reconstructions are obtained after direct application AGCM.}}}
\label{fig:7}
\end{figure*}

% Figure \ref{fig:fig3} presents typical behavior of noisy
% backscattered data for our two objects.
 Figures \ref{fig:4}, \ref{fig:5} present reconstructions of
 $\tilde{c}$ which we have obtained on three time adaptively refined
 mesh inside the domain $G_1$ for different values of $\omega$ in
 (\ref{f}).  We note that different initial guesses $\tilde{c}_0$ in
 \eqref{functional} produce different symmetric structures inside
 $G_1$ with different values of the function $\tilde{c}(x)$, compare
 reconstructions presented on Figures \ref{fig:7}.  Left images of Figure
 \ref{fig:7} present reconstructions obtained in ACGM
 when the optimized function $\tilde{c}$, obtained on a coarse mesh, is
 sequentially interpolated on the one, two and three times refined
 mesh. Then this interpolated function is taken as an initial guess
 $\tilde{c}_0$ in \eqref{functional} and optimized further to get
 reconstruction on the third refined mesh.  Right images of Figure
 \ref{fig:7} are obtained after direct application of the  adaptive
 algorithm of Section \ref{sec:ad_alg}.  Optimized values of
 $\tilde{c}(x)$ obtained on Figures \ref{fig:4}--\ref{fig:7} can be of
 physical interest since they present symmetric structured domains
 with almost the same material in every structured layer.

\section*{Acknowledgments}

%This research is supported by the founding from Area of Advance
%``Nanoscience and Nanotechnology'', Chalmers University of Technology.
  The research of L.B. is
  supported by the sabbatical programme at the Faculty of Science,
  University of Gothenburg. The research of E.S. is  supported by the Ministry of Education and Science of the Russian Federation, Project No.~1.894.2017/$\Pi$.

\emergencystretch=\hsize

\begin{center}
\rule{6 cm}{0.02 cm}
\end{center}

%\end{floatingfigure}

\begin{thebibliography}{99}


 \bibitem{BKS}
 \newblock Bakushinsky A., Kokurin M.Y.,  Smirnova A.,
 \newblock \emph{Iterative Methods for Ill-posed Problems},
 \newblock Inverse and Ill-Posed Problems Series 54, De Gruyter, 2011.


 \bibitem{B}   L.~Beilina, Adaptive hybrid FEM/FDM methods for inverse
  scattering problems.  \emph{ Inverse Problems and Information
    Technologies}, V.1, N.3, 73-116, 2002.


\bibitem{hybrid} L. Beilina, Domain decomposition finite
  element/finite difference method for the conductivity reconstruction
  in a hyperbolic equation, \emph{Communications in Nonlinear Science
    and Numerical Simulation}, Elsevier, 37, p.222-237, 2016.

% \bibitem{BMaxwell}   L. Beilina,
%  Energy estimates and numerical  verification of the
%   stabilized domain decomposition finite  element/finite difference
%   approach for time-dependent Maxwell's  system,
% \emph{Cent. Eur. J. Math.},  11, 702-733, 2013.


\bibitem{BJ} L. Beilina and C. Johnson, A posteriori error estimation in
computational inverse scattering, \emph{Mathematical Models in Applied
Sciences}, 1, 23-35, 2005.

\bibitem{BCL} L. Beilina, M. Cristofol, S. Li, Uniqueness and
  stability of time and space-dependent conductivity in a hyperbolic
  cylindrical domain, arXiv:1607.01615.

%\bibitem{BCN}  L. Beilina, M. Cristofol and K. Niinim\"aki, Optimization
%  approach for the simultaneous reconstruction of the dielectric
%  permittivity and magnetic permeability functions from limited
%  observations, \emph{Inverse Problems and Imaging}, 9 (1), pp. 1-25, 2015.

%\bibitem{BN} L. Beilina and K. Niinim\"aki, Numerical studies of the
%  Lagrangian approach for reconstruction of the conductivity in a
%  waveguide,  arXiv:1510.00499, 2015.


\bibitem{btkm14} L.~Beilina, Nguyen~T.T., M.~Klibanov, and
  J.~Malmberg, Reconstruction of shapes and refractive indices from
  backscattering experimental data using the adaptivity, \emph{Inverse
    Problems} \textbf{30}, 105007 2014.

%\bibitem{btkm14b}
%L.~Beilina, Nguyen~T.T., M.~Klibanov, and J.~Malmberg,   Globally convergent an%d adaptive finite element
%  methods in imaging of buried objects from experimental
%  backscattering radar measurements,\emph{J. Comput. Appl.
%  Math.}, 289, pp. 371-301, 2015, doi:10.1016/j.cam.2014.11.055.



  \bibitem{BOOK}  L. Beilina, M.V. Klibanov,
    \emph{Approximate global convergence and adaptivity for coefficient inverse problems},
    Springer, New-York, 2012.

\bibitem{AoA} L. Beilina, L. Mpinganzima, P. Tassin, Adaptive
  optimization algorithm for the computational design of nanophotonic
  structures, \emph{IEEE, Proceedings of the 2016 International
    Conference on Electromagnetics in Advanced Applications, ICEAA
    2016}, pp. 420-423, 2016, doi:10.1109/ICEAA.2016.7731416.


\bibitem{Brenner} S.~C.~Brenner and L.~R.~Scott, \emph{The
  Mathematical Theory of Finite Element Methods}, Springer-Verlag,
  Berlin, 1994.

%\bibitem{Cohen} G. C. Cohen, \emph{Higher  Order Numerical Methods for
%  Transient Wave Equations}, Springer-Verlag, Berlin, 2002.


\bibitem{EM}  B. Engquist  and A. Majda, Absorbing boundary conditions for
the numerical simulation of waves, \emph{\ Math. Comp.}, 31, 629-651, 1977.


\bibitem{Engl}  H. W. Engl,  M. Hanke and A. Neubauer,  \emph{\ Regularization
of Inverse Problems}, Kluwer Academic Publishers, Boston, 2000.

%\bibitem{JoanJohnson} Joannopoulos, Johnson, Winn and Meade,
%  \emph{Photonic Crystals: Molding the Flow of Light}, Second edition,
%  Princeton Univ. Press, 2008.




%\bibitem{KBKSNF} Kuzhuget, A.V., Beilina, L., Klibanov, M.V.,
%  Sullivan, A., Nguyen, L., Fiddy, M.A.,  Blind experimental data
%  collected in the field and an approximately globally convergent
%  inverse algorithm, \emph{Inverse Problems}, V.28, N.9, 2012,
%  DOI:10.1088/0266-5611/28/9/095007

\bibitem{lad} O. A. Ladyzhenskaya, \emph{Boundary Value Problems of
  Mathematical Physics}, Springer-Verlag, Berlin, 1985.

% \bibitem{Maier} Maier,\emph{ Plasmonics: Fundamentals and Applications}, Springer, 2007.

\bibitem{cloak1} U. Leonhardt, \emph{Optical Conformal
  Mapping}, Science, 312, pp. 1777-1780, 2006.

\bibitem{cloak2} J. B. Pendry, D. Schurig and D. R. Smith, \emph{Controlling
  electromagnetic fields}, Science, 312,pp. 1780 - 1782, 2006.

\bibitem{cloak3} H. Chen, C. T. Chan, \emph{Acoustic cloaking and
  transformational acoustics}, Journal of Physics, IOP
  Publishing,  2010, doi:10.1088/0022-3727/43/11/113001.

 \bibitem{petsc}
 \newblock PETSc, Portable, Extensible Toolkit for
 Scientific Computation, http://www.mcs.anl.gov/petsc/

\bibitem{Peron} O.Pironneau, \emph{Optimal Shape Design for Elliptic
  Systems}, Springer-Verlag, Berlin, 1984.


%\bibitem{SoukoulisW}  Soukoulis, & Wegener, Nature Photon. 5, 523, 2011.


%\bibitem{NBKF}  N. T. Th\`{a}nh, L.~Beilina, M.~V. Klibanov and M.~A.
%Fiddy, Reconstruction of the refractive index from experimental
%backscattering data using a globally convergent inverse method, \emph{SIAM
%J. Scientific Computing}, 36 (3), pp.273-293, 2014.

%\bibitem{NBKF2}  N. T. Th\`{a}nh, L. Beilina, M. V. Klibanov, M. A. Fiddy,
%  Imaging of buried objects from experimental backscattering
%  time-dependent measurements using a globally convergent inverse
%  algorithm, \emph{SIAM Journal on Imaging Sciences}, 8(1), 757-786,
%  2015.


\bibitem{tikhonov} A. N. Tikhonov, A. V. Goncharsky, V. V. Stepanov  and
 A. G.  Yagola, \emph{Numerical Methods for the Solution of
    Ill-Posed Problems},  Kluwer, London, 1995.



\bibitem{waves}  WavES, the software package, http://www.waves24.com~.

%\bibitem{ZheludevK} Zheludev, Kivshar, Nature Mater. 11, 917, 2012.

\bibitem{CFL67}{R. Courant, K. Friedrichs and H. Lewy}, On the partial
  differential equations od mathematical physics, \emph{IBM Journal of
    Research and Development}, 11(2),  215-234, 1967.


%%%%%%%%%%%%%%%%%%%%%%%%%%%%%%%%%%%%%%



\end{thebibliography}
\end{document}